\documentclass[10pt]{article}

\usepackage[dvips]{graphicx}
\usepackage{amsmath,amssymb,epsfig}

\title{Some extensions of the class of $k$-convex bodies}
\author{V. Golubyatnikov\footnote{ Supported by
President of Russia grant for leading scientific schools No.~8526.2006.1
and }
\\
 {\em\small Sobolev Institute of Mathematics, Novosibirsk, 630090, Russia} \\
 {\em\small e-mail: glbtn@math.\,nsc.\,ru}\\
 V. Rovenski\footnote{grant FP7-PEOPLE-2007-2-1-IEF, No. 219696 of Marie-Curie action.}
  {\em\small
 University of Haifa, Haifa, 31905, Israel}\\
 {\em\small e-mail: rovenski@math.\,haifa.\,ac.\,il}
 }

\date{}

\begin{document}

\newtheorem{theo}{Theorem}
\newtheorem{defi}{Definition}
\newtheorem{lem}{Lemma}
\newtheorem{prob}{Problem}
\newtheorem{prop}{Proposition}
\newtheorem{cor}{Corollary}
\newtheorem{exam}{Example}
\newtheorem{rem}{Remark}
\newtheorem{exer}{Exercise}

\maketitle


\textbf{Abstract}.\footnote{\textit{Keywords and phrases}: $k$-convex, $k$-visible body, supporting ball, circular projections, visual hull, geometric tomography

\textit{AMS Mathematics Subject Classification}: 52A01, 52A30, 52A35}
We study relations of some classes of $k$-convex, $k$-visible bodies in Euclidean spaces. We introduce and study \textrm{circular projections} in normed linear spaces and classes of bodies related with families of such maps, in particular, \textrm{$k$-circular convex} and \textrm{$k$-circular visible} ones.
 Investigation of these bodies more general than $k$-convex and $k$-visible ones allows us to generalize some classical results of geometric tomography and find their new applications.

\section*{Introduction}

As usual, a convex body in $n$-dimensional Euclidean space ${\Bbb{E}}^n$ (over reals, complex or quaternion numbers) is a compact convex $K$ set with non-empty interior, $int(K)$, whose closure coincides with the body, $cl(int(K))=K$.
The classical \textbf{problem} of geometric tomography concerns the reconstructing a body $K$ of ${\Bbb{E}}^n$ from the knowledge of its orthogonal projections onto $(n-k)$-planes of a preselected family~$\mathcal{P}_k$.

Convex bodies are characterized by the property that they can be represented as intersections of families fo (supporting) half-spaces.
The axiomatic approach to the notion of convexity consists of the following. In an arbitrary set $X$, we choose a certain family $\mathcal{B}$ of subsets (called the convexity base).
A~set $K$ is called $\mathcal{B}$-convex if it is an intersection of some subfamily of sets from $\mathcal{B}$.

Reshetnyak \cite{resh56} considered some generalization of convex surfaces, called $\delta$-touch\-ed surfaces, in order to study behavior of geodesic lines on these surfaces.
Following the same direction we have introduced in \cite{vr} the $\varepsilon$-convex bodies (more general than convex ones) and investigated the problem of determination of $\varepsilon$-convex bodies by their projection-type images.
Instead of supporting half-spaces (as a  family $\mathcal{B}$) we used there some different shapes such as complements to balls and cylinders and introduced \textrm{circular projections} onto punctured hyperplanes.

Some results for the convex bodies can be proved for $k$-convex and $k$-visible bodies (see \cite{vgol}) and $k$-visual bodies (see \cite{MU67}, \cite{LM70}) as well.
The definition of a $k$-convex submanifold was introduced in
\cite{Shefel'_1969} as a multi-dimensional generalization of a convex surface.
In order to prove stability of determination of $k$-visible bodies by their images under projections, the first author \cite{vgol92} introduced $\alpha$-rough bodies, where $\alpha\in(0,\pi]$; for $\alpha=\pi$ we get the class of convex bodies.

In the present work we extend the classes of $k$-convex and $k$-visible bodies in ${\Bbb{E}}^n$
and study the problems of their reconstruction from the knowledge of orthogonal projections onto $(n-k)$-planes.
We introduce our main tool, the families $\mathcal{P}_{k}$ of \textbf{circular projections} onto punctured $(n-k)$-dimensional planes.

In order to extend our previous results \cite{vr}, \cite{vgol}, \cite{vgol99}, \cite{Rov_2006}, following to the "soft-hard" ranking of geometrical categories described in~\cite{gro}, we study the classes of $\mathcal{C}$-\textbf{convex} $\mathcal{KP}^{k}_i,\mathcal{K}^{k,\varepsilon}_i$
and $\mathcal{VP}^{k}_i,\mathcal{V}^{k,\varepsilon}_i$ ($i=1,2,\ \varepsilon>0$) of $\mathcal{C}$-\textbf{visible} bodies,
more general than $\varepsilon$-convex or $\varepsilon$-visible classes of bodies introduced in \cite{vr}.
 In~the definitions of these classes we consider

 -- the boundary points of the bodies ($i=1$),

 -- the points disjoint from the bodies ($i=2$).

\noindent
We prove strong inclusions between these classes. Note that the classes $\mathcal{K}^{k}_2$, $\mathcal{KP}^{k}_2$ and $\mathcal{K}^{k,\,\varepsilon}_2$ are closed under intersections of bodies,
some results for the second type of bodies ($i=2$) can be extended to the wider 1st one ($i=1$).

It is interesting to find conditions \textbf{C} such that
"\textit{If the bodies $U,V{\subset}{\Bbb{E}}^n$ satisfy \textbf{C} and their circular projections of a preselected family $\mathcal{P}_k$ coincide then $U=V$}".

Much more hard and general \textbf{question} is the following:
What kind of condition \textbf{CT} and what kind of transformations \textbf{PT, ST} should be taken in order to get a proposition:
"\textit{If $ U, V \subset{\Bbb{E}}^n $ satisfy
 \textbf{CT} and their circular projections
of a family $\mathcal{P}_k$ are \textbf{PT}-equivalent, then the bodies $U, V$ are \textbf{ST}-equivalent in the ambient space}"?

For the case of orthogonal projections some answers to the question above were obtained in
\cite{vgol}, \cite{vgol99}.
The study of circular projections and circular $\mathcal{C}$-visible bodies
allows us to extend some classical results of the theory of convex bodies, see for example  \cite{vgol}, \cite{rgar06}, and to find their new applications in geometric tomography.

\section{$k$-Convex bodies and visual hulls of sets}

In what follows we as usual consider compact bodies in ${\Bbb{E}}^n$ over reals.
Let $P^k\subset{\Bbb{E}}^n$ be a $k$-dimensional plane (i.e., an affine $k$-subspace). Denote~by
\begin{equation}\label{E-tube}
 B(P^k,r)=\{x\in{\Bbb{E}}^n: {\rm dist}(x,P^k)\le r\}
\end{equation}
a {solid cylinder of radius $r>0$ with the axis $P^k$}, that is the metric product of ${\Bbb{E}}^k$ and a ball $B^{n-k}(O,r)\subset{\Bbb{E}}^{n-k}$. For $k=0$, (\ref{E-tube}) gives us a ball $B^{n}(C,r)$.

\begin{defi}\label{D-supp}\rm
Let $L$ be a body (usually we consider balls or cylinders) or a $k$-dimensional plane in ${\Bbb{E}}^n$. We~call $L$ \textbf{supporting} a body $K$ if $L$ intersects $\partial K$ and $L$ is disjoint from int\,$K$.

A half-plain $P^{k+1}_{^\ge}$ can be represented as the union $P^{k+1}_{^>}\cup P^{k}$ of an open half-plane $P^{k+1}_{^>}$ and a boundary $k$-plane.
A~half-plane $P^{k+1}_{^\ge}$ will be called \textbf{supporting} a body $K$ if $P^{k}\cap\partial K\ne\emptyset$ and $P^{k+1}_{^>}\cap B(P^{k},r)$ is disjoint from int\,$K$ for some $r>0$ depended only on $K$.
In this case, obviously, $P^{k}$ is supporting~$K$.
In~what follows  for simplicity we assume $r=\infty$.
\end{defi}

\begin{defi}\label{D-Kki}\rm
A body $K\subset{\Bbb{E}}^n\ (n>2)$ is called \textbf{$k$-convex of a class $\mathcal{K}^{k}_i$} if

\hskip-1mm
$\mathcal{K}^{k}_1$: \textit{any point $x\in\partial K$ belongs to
the boundary of a half-plane $P^{k+1}_{^\ge}$ supporting~$K$}.

\hskip-1mm
$\mathcal{K}^{k}_2$: \textit{any point $x\notin K$ belongs to
the boundary of a half-plane $P^{k+1}_{^\ge}$ disjoint from $K$}.

\noindent
A body $K\subset{\Bbb{E}}^n$ is called $k$-\textbf{visible of a class $\mathcal{V}^{k}_i$} if

\hskip-1mm
$\mathcal{V}^{k}_1$: \textit{any plane $Q^{k-m}$ supporting $K$ belongs to
the boundary of a half-plane $P^{k+1}_{^\ge}$ supporting~$K$}.

\hskip-1mm
$\mathcal{V}^{k}_2$: \textit{any plane $Q^{k-m}$, which is disjoint from $K$, belongs to
the boundary of a half-plane $P^{k+1}_{^\ge}$ disjoint from $K$}.

\vskip.1mm\noindent
Here $0<m\le k<n$, and the $0$-dimensional plane is a point.
\vskip.5mm\noindent
A~connected boundary component of a $k$-convex (or $k$-visible) body $K$ is called a $k$-\textbf{convex}
(or $k$-\textbf{visible}) \textbf{hypersurface} of a certain class listed above.
\end{defi}

\begin{rem}\label{R-Kki1}\rm
(a) A $\Pi$-shaped non-convex body in ${\Bbb{E}}^3$ belongs to $\mathcal{K}^{k}_1$.

(b) A connected body of a class $\mathcal{K}^{n-1}_1$ is a convex~one.
\end{rem}

\begin{prop}\label{R-Kki}
The conditions $\mathcal{K}^{k}_i\ (i=1,2)$ are equivalent to the next ones:

\hskip-3mm $(i=1)$ \textit{for any $x\in\partial K$,
 $K\in{\mathcal K}^k_1 $, there is orthogonal projection $f:{\Bbb{E}}^n\to P^{n-k}$
 such that $f(x)$ is a vertex of a ray supporting $f(K)$}.

\hskip-3mm $(i=2)$ \textit{for any $x\notin K$, $K\in{\mathcal K}^k_2 $, there is
orthogonal projection $f:{\Bbb{E}}^n\to P^{n-k}$ such that $f(x)$
 is a vertex of a ray disjoint from $f(K)$}.

\noindent
The condition $\mathcal{V}^{k}_i\ (i=1,2)$ is equivalent to the following one:

\hskip-3mm
$(i=1)$ \textit{for any plane $Q^{k-m}$ supporting $K$, $K\in\mathcal{V}^{k}_1$, there is
orthogonal projection $f:{\Bbb{E}}^n{\to} P^{n-k}$ such that $f(Q^{k-m})$ is a vertex of a ray supporting $f(K)$}.

\hskip-3.5mm
$(i=2)$ \textit{for any plane $Q^{k-m}$ disjoint from $K$, $K{\in}\mathcal{V}^{k}_2$, there is
orthogonal projection $f:{\Bbb{E}}^n{\to} P^{n-k}$ such that
$f(Q^{k-m})$ is a vertex of a ray disjoint from~$f(K)$}.
\end{prop}

The proof of Proposition~\ref{R-Kki} follows from the definitions.

\begin{rem}\rm
The classes of $k$-\textbf{convex} and $k$-\textbf{visible} bodies can be defined relative to any non-empty family $\mathcal{P}_{k}$ of orthogonal projections onto $(n-k)$-planes.
\end{rem}

Class $\mathcal{V}^{k}_2$ for $k=1$ is closed under the intersection of bodies in ${\Bbb{E}}^n$, for $k>1$ this claim fails, \cite{vgol}.
Note that $\mathcal{V}^{1}_i=\mathcal{K}^{1}_i\ (i=1,2)$.
 Remark that the properties $(i=1)$
fail if one will replace the definitions of $\mathcal{K}^{k}_1$ and $\mathcal{V}^{k}_1$ by weaker ones:

\hskip-3.5mm
$\mathcal{K}^{k}_{1-}$: \textit{any point $x\in\partial K$ belongs to a $k$-plane $P^{k}$ supporting $K$}.

\hskip-3.5mm
$\mathcal{V}^{k}_{1-}$: \textit{any plane $Q^{k-m}$ supporting $K$ belongs to a plane $P^{k}$ supporting $K$ as~well}.
\newline
Similarly one may formulate weaker conditions $\mathcal{K}^{k}_{2-}$ and $\mathcal{V}^{k}_{2-}$.
Some of results of this section can be proved for the above classes of bodies.

\begin{exam}\label{Ex-1-2A}\rm
\nopagebreak
 (a) The union of black squares of a chess-board and its $n$-dimensio\-nal analogues provide examples of bodies  (not homeomorphic to a ball) in ${\Bbb{E}}^n$ of the classes $\mathcal{K}^{n-1}_{1-}\setminus\mathcal{K}^{n-1}_1$
 and $\mathcal{V}^{n-1}_{1-}\setminus\mathcal{V}^{n-1}_1$.
 Example with homeomorphic to a ball bodies is given in what follows.

 (b) Consider a homeomorphic to a ball body $K$ in ${\Bbb{E}}^3$ between two helicoids,
 surfaces $M_1$ and $M_0$, where $M_h: [u\cos v, u\sin v, v+h]$ and $0\le v\le2\pi, 0\le u\le1$.
The projection of $K$ onto $xy$-plane is a unit disc.
The property ($\mathcal{K}^{k}_{1-}$) for $k=1$ is satisfied, but the $xy$-projection of the $z$-axis
(supporting $K$) is the inner point (the origin) of the disc.
Hence, $K\not\in\mathcal{K}^{k}_{1}$.
\end{exam}

The following result was proved in \cite{vgol} for a smaller class $\mathcal{V}^{k}_2$.

\begin{lem}
(a) A connected body $K\in\mathcal{V}^{n-1}_1$ in ${\Bbb{E}}^n$ is convex.
(b) The projection of a body $K\in\mathcal{V}^{k}_1$ in ${\Bbb{E}}^n$ onto any hyperplane
    belongs to a class $\mathcal{V}^{k-1}_1$. The corresponding visibility will be for its projections onto planes of smaller dimension, in particular, the projection of a body $K\in\mathcal{V}^{n-2}_1$ in ${\Bbb{E}}^n$ onto any $3$-plane belongs to a class $\mathcal{V}^{1}_1$.
\end{lem}

\textbf{Proof}.
 Part (a) of lemma follows from the fact that the body $K$ lies on one side of its supporting hyperplane.

(b) Let $K(\omega)$ be the projection of $K$ onto a hyperplane $P^{n-1}(\omega)$, orthogonal to unit vector $\omega\in{\Bbb{E}}^n$ and let $y\in\partial K(\omega)$. All inverse images of this point under projecting lie on the boundary of $K$. Let $z\in\partial K$ be one of inverse images of the point $y$. The line $l(yz)$ is supporting $K$ and, by definition of the class $\mathcal{V}^{k}_1$, belongs to $(k+1)$-dimensional half-plane $H^{k+1}$, supporting $K$.
The image of $H^{k+1}$ under projection onto $P^{n-1}(\omega)$ is a $k$-dimensional half-plane, supporting $K(\omega)$.$\,\square$

\begin{theo}\label{P-separate33}
 $
 \mathcal{K}^{k}_2\subset\mathcal{K}^{k}_1$, \ \
 $
 \mathcal{V}^{k}_2\subset\mathcal{V}^{k}_1$
 \ and \
 $\mathcal{V}^{k}_i\subset\mathcal{K}^{k}_i\ (i=1,2)$.
\end{theo}

\textbf{Proof}.
$\mathcal{K}^k_2\subset\mathcal{K}^k_1$.
Given $K\in\mathcal{K}^k_2$ and $\bar x\in\partial K$, consider a sequence $x_i\not\in K$ such that $x_i\to\bar x$ for $i\to\infty$. Due to $(\mathcal{K}^k_2)$, for each $i$ there is a $k$-plane $P^{k}_i$ containing $x_i$ which
bounds a half-plane $P^{k+1}_{^\ge,i}$ disjoint from $K$.
Let $p_i$ be the normal to $P^{k}_i$ at $x_i$ directed inside of $P^{k+1}_{^\ge,i}$.
Taking a subsequence of integers, one may assume that there are a unit vector $\bar p=\lim\limits_{i\to\infty}p_i$ and a $k$-plane $\bar P^k=\lim\limits_{i\to\infty}P^k_i$ containing $\bar x$. Then $\bar P^k$ bounds a half-plane $\bar P^{k+1}_{^>}$
and $\bar p\,\perp\,\bar P^k$ is a unit normal to $\bar P^{k}$ at $\bar x$ directed inside of $\bar P^{k+1}_{^\ge}$. We claim that $\bar P^{k+1}_{^\ge}$ is supporting $K$.
To show this assume an opposite, that is $P^{k+1}_{^>}$
contains a point $y\in{\rm int}\,K$.
Using standard isometry $\phi_i: \bar P^{k+1}_{^\ge}\to P^{k+1}_{^\ge,i}$,
$\phi_i(\bar x)=x_i$,
we find a sequence of points $y_i=\phi_i(y)\notin\,K$ converging to $y$.
From this it follows that $y\notin{\rm int}\,K$, that is a contradiction.
 Other inclusions can be proved analogously.
Strong inclusions follow from examples in the sequel.$\,\square$

\begin{exam}\label{Ex-1-2}\rm
\nopagebreak

(a) $K\in\mathcal{K}^{k}_1\setminus\mathcal{K}^{k}_2$ for $k=1,\,n=3$.
 Consider an envelope $D$ (homeomorphic to a disc) with four unit square faces, and two sections
 along $p$ and $q$, see Fig.~\ref{ex-K12}. Let us bend $D\times[-c,c]$ (for small $c$) about $r$ by an angle a bit smaller than $90^o$, and by $90^o$ about the lines $p,q$.
 We obtain a body $K\subset{\Bbb{E}}^3$ that partially surrounds a cube (from four faces).
 Obviously, $K\in\mathcal{K}^{1}_1$, but the property $\mathcal{K}^{1}_2$ is not satisfied
 for the point $C(0,-1/4,1/4)\notin K$ close to the front face~($r$).
 A~cylinder $K'=K\times[-a,a]\subset{\Bbb{E}}^4$ is the body $K'\in\mathcal{V}^{k}_1\setminus\mathcal{V}^{k}_2$ for $k=2,\,n=4$.

(b) Example of a body $K\in\mathcal{K}^{k}_2\setminus\mathcal{V}^{k}_2$ for $k=2,\,n=4$,
(i.e., a 2-convex body in ${\Bbb{E}}^4$ that is not 2-visible) is given in \cite{vgol}.
\begin{figure}
\begin{center}
\includegraphics[scale= 0.3,angle=0]{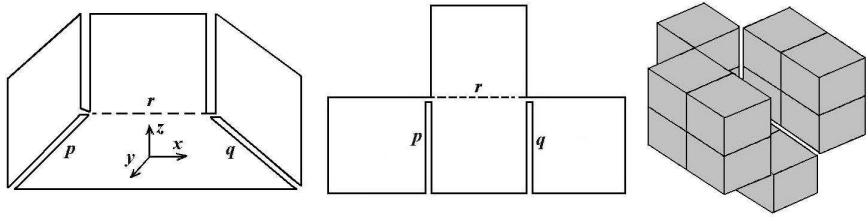}
\caption{\small A body $K\in\mathcal{K}^{k}_1\setminus\mathcal{K}^{k}_2\ (k=1)$, its "envelope" and model of 13 cubes.}
\label{ex-K12}
\end{center}
\end{figure}
\end{exam}

The following simple result was obtained by the first author on the initial stage of the study of the question that was formulated in Introduction.

\begin{lem}[\cite{vgol}]\label{L-ortG}
Let $V_1, V_2\in{\cal K}^{k}_2$ be compact bodies in ${\Bbb{E}}^n, n>2,\ n-k>1$ (real, complex or quaternion)
and $G$ be a group of all either translations, or homotheties or is trivial.
If the projections of $V_1, V_2$ onto any $(n-k)$-plane $P$ coincide relative to $G$, then $V_1,V_2$ also coincide relative to $G$.
\end{lem}

\hskip-1mm
This claim is wrong in the real case for $n-k=1$ and non-trivial $G$, since two different bodies of constant width have isometric projections onto any line.

\begin{rem}\label{R-03}\rm
Let $K\in\mathcal{K}^k_1\ (k>n/2)$ be a body in ${\Bbb{E}}^n$, and let $\partial K$ be a $C^2$-regular hypersurface with a unit normal $\xi$ directed inside.
By conditions, any point $x\in\partial K$ belongs to a plane $P^k_x$ disjoint from int\,$K$.
From this we conclude that $P^k_x\subset (TM)_x$,
and the 2nd fundamental form of $\partial K$ is non-negative on~$P^k_x$.
hypersurface. Hence at least $k$ eigenvalues of the (symmetric) Weingarten operator $A_\xi$
are non-negative. In this case, the homology groups (with integer coefficients) $H_{n-1-k+s}(\partial K)=0$ for $0<s<2k-n-1$,
see Lemma 16 in \cite{Bor82} and also Lemma 6(b) in \cite{Rov_2006}.
\end{rem}

 \begin{defi}\label{D-visualh}\rm
Given set $W$ in ${\Bbb{E}}^n$ and a non-empty family $\mathcal{P}_{k}$ (of orthogonal projections onto $(n-k)$-planes), a \textbf{visual hull} $\mathcal{P}_{k}\langle W\rangle$ is the largest set~$V$ such that $f(V)\subseteq f(W)$ for all $f\in\mathcal{P}_{k}$.
 A~body $K\subset{\Bbb{E}}^n$
is \textbf{visual} relative to $\mathcal{P}_{k}$ if $\mathcal{P}_{k}\langle K\rangle=K$~\cite{MU67}.
\end{defi}

If $f$ is a projection onto subspace $P\subset{\Bbb{E}}^n$, $\mathcal{R}(f)$ denotes its range, i.e., $P$.
Denote by $\bar{\mathcal{P}}_{k}$ a family of orthogonal projections onto all $(n-k)$-planes in ${\Bbb{E}}^n$.

\begin{prop}[\cite{MU67}]
$\,\mathcal{P}_{k}\langle K\rangle=K$ if and only if for every $x\notin K$
there exists a $f\in\mathcal{P}_k$ such that $x+\mathcal{R}(f)^\perp$ does not intersect~$K$.
\end{prop}

The \textbf{visual hulls} of bodies in ${\Bbb{E}}^3$ for finite families $\mathcal{P}_k$
(reconstructions of 3-D shapes from multiple calibrated photographs) are studied
in computational geometry and computer vision, see \cite{LM70}, \cite{LF07}, etc.

\begin{prop}\label{P-003} $K\in\mathcal{K}^{k}_2$ if and only if $\,\bar{\mathcal{P}}_{n-k}\langle K\rangle=K$.
\end{prop}

\textbf{Proof}. Suppose that there is $x\in\bar{\mathcal{P}}_{n-k}\langle K\rangle\setminus K$. By definition of $\mathcal{K}^{k}_2$, $x$ belongs to some $k$-plane $P^k$ which is disjoint from $K$, hence the projection of $K$ onto orthogonal $(n-k)$-plane $(P^{k})^\perp$ does not contain the image of $x$.
Thus $x\notin\bar{\mathcal{P}}_{n-k}\langle K\rangle$, a contradiction.$\,\square$

\section{Circular projections and convexity}

We will extend the classes of $k$-convex bodies in ${\Bbb{E}}^n$ in a way that the role of $k$-planes will play the $k$-dimensional spheres of radius $1/\varepsilon$ for a given $\varepsilon>0$, and the circular projections will be used instead of orthogonal ones.

\subsection{Circular projections}

A $(n-k)$-plane $P(C)$ containing a point $C\in{\Bbb{E}}^n$ is called the
\textbf{punctured $(n-k)$-plane}. Denote by $P^{\perp}(C)$ the
punctured $k$-plane orthogonal to $P(C)$. Let $S^{n-1}(C,r)$ be a
sphere of radius $r>0$ centered in $C$, and $S^{n-k-1}(P(C),r)=S^{n-1}(C,r)\cap P(C)$ be a sphere of the same radius and center.
 Given $x\notin P^{\perp}(C)$, let $S^\perp(P(C),x)$ be a large
 $k$-sphere in a sphere $S^{n-1}(C,|Cx|)$ through $x$ and orthogonal
 to a large sphere $S^{n-k-1}(P(C),|Cx|)$.

\begin{defi}\label{D-circ-pro2}\rm
 Let $P(C)$ be a punctured $(n-k)$-plane.
 Define a \textbf{circular projection} $f_{P(C)}:{\Bbb{E}}^n{\setminus} P^{\perp}(C)\to P(C){\setminus}C$ as follows. Given $x\notin P^{\perp}(C)$, let $f_{P(C)}(x)$ be the closest to $x$ point of the sphere $S^{n-k-1}(P(C),|Cx|) \subset P(C)$.
 \end{defi}

 One may represent a circular projection explicitly by a formula.
If $C=O$, $\omega = {\rm pr}_{P^\perp(C)}(x)/||\,{\rm pr}_{P^\perp(C)}(x)||$ and $x\not\in\,P^\perp(O)$ then
\begin{equation}\label{E-fomega}
\begin{array}{c}
 f_{P(O)}(x)=\frac{\|x\|}{\sqrt{x^2-(\omega,x)^2}}[x-(\omega,x)\,\omega].
\end{array}
\end{equation}

 The circular projection $f_{P(C)}$ keeps the distance of a point $x$ to $C$.
 The~circular projection of a ball can be a non-convex body, see \cite{vr}.

The circular projections onto punctured 1-planes (lines) have more simple 'nature' than ones onto higher dimensional planes; circular projection onto any two lines through $C$ are similar in some sense. Here, one can define a \textbf{positive circular projection}
$|f_C|:{\Bbb{E}}^n\to{\mathbb R}_+ $ as $|f_C|(x):=|Cx|$.

For the real Euclidean space ${\Bbb{E}}^n$, the points $|f_C|(x)$ and $x$ belong to the same half-space with respect to hyperplane $P^{\bot}(C)$. For the complex or the quaternion space ${\Bbb{E}}^n$ the intersection (with a line) $P(C)\cap S^{an-1}(C,|Cx|)\ (a=2,4$) is a $(a-1)$-sphere, and $|f_C|(x)$ is the nearest to $x$ point of this sphere.

\begin{rem}\rm
In order to generalize the notion of circular projections we consider fibration
$\pi_0:{\Bbb{E}}^n\setminus(0\times{\Bbb{E}}^{k})\to{\Bbb{E}}^{n-k}\setminus\{0\}$, the fiber $S_0(x)=\pi^{-1}(x)$ is a $k$-dimensional surface,
in particular case, such construction has ellipsoidal fibers.

Let $g:[-a,b]\to{\Bbb{R}}_+$ (where $0<a,b\le\infty$) be a smooth convex function such that $g(0)=1,\,g'(0)=0$ and the curvature $k_g\le1$. Denote by $\gamma_R$ a curve $x=R\,g(y/R)$ that is either contained in the interior of a real half-space ${\Bbb{E}}^2_+$ or meets its boundary $m$ transversally.

Suppose that ${\Bbb{E}}^n={\Bbb{E}}^{n-k}\times{\Bbb{E}}^{k}$ for some $k$, and let ${\Bbb{E}}^2_+$ is orthogonal to both factors
(say, $x$-axis belongs to ${\Bbb{E}}^{n-k}$ and $y$-axis belongs to ${\Bbb{E}}^{k}$).
Now, the group $SO(n)$ contains a subgroup $G$ isomorphic to $SO(n{-}k){\times} SO(k)$ whose action on ${\Bbb{E}}^n$ splits
along the factors $G_1=SO(n{-}k)$ and $G_2=SO(k)$.
This gives rise to a $G$-invariant hypersurface $M_R^{n-1}=G(\gamma_R)\subset{\Bbb{E}}^n$ called a
\textbf{hypersurface of revolution with $\gamma_R$ as profile}. Each $M_R^{n-1}$ intersects ${\Bbb{E}}^{n-k}\times\{0\}$
by a sphere $G_1(\gamma_R(0))=S^{n-k-1}(R)$ of radius $R$, i.e., the $G_1$-orbit of a point $\gamma_R(0)$.
The action of a subgroup $G_2$ gives us a fibration $\pi_0:{\Bbb{E}}^n\setminus (0\times{\Bbb{E}}^{k})\to {\Bbb{E}}^{n-k}\setminus\{0\}$,
the fiber $S_0=G_2(\gamma_R)$ is the $k$-dimensional surface. Moreover, $M_R^{n-1}$ is a union of such surfaces.

For example, if $g(y)=\sqrt{1-x^2}$ then $\gamma_R$ is a semi-circle and $M_R^{n-1}$ is a sphere of radius $R$,
and $G_1(\gamma_R)$ is a large $k$-sphere orthogonal to $(n-k)$-sphere $M_R^{n-1}\cap{\Bbb{E}}^{n-k}$.
Similarly, $g(y)=(b/a)\sqrt{1-x^2}$ produces the ellipsoids.
\end{rem}

\begin{defi}\label{D-Gcirc-pro2}\rm
Given a punctured $(n-k)$-plane $P(C)$, consider a Euclidean motion
 $T:{\Bbb{E}}^n{\setminus} P^{\perp}(C)\to{\Bbb{E}}^{n}\setminus{\Bbb{E}}^k$.
Let $g:[-a,b]\to{\Bbb{R}}_+$ be as above.
Then $\pi_{C}:=T^{-1}\circ\pi_0\circ T:{\Bbb{E}}^n{\setminus}P^{\perp}(C)\to P(C)\setminus C$ is a fibration by $k$-dimensional surfaces. Denote by $M^{k}(P(C),x)$ such a fiber through a point $x\not\in P^\perp(C)$. For any such point denote by $f_{P(C)}(x)$ the closest to $x$ point of the sphere $S(P(C),g,x):=S(C,R(x))\cap P(C)$, where $M^{n-1}_{C,R(x)}$ a hypersurface containing $x$.
 We call $f_{P(C)}:{\Bbb{E}}^n{\setminus} P^{\perp}(C)\to P(C){\setminus} C$
 a $g$-\textbf{non-linear projection} onto $P(C)$. Denote by $S^\perp(P(C),x)$
 a $k$-dimensional 'surface-meridian' in $M^{n-1}_{C,R(x)}$ through $x$
 (orthogonal to a sphere $S(P(C),g,x)$).
 \end{defi}

The definitions (starting form Definition~\ref{D-CVH}) and results in what follows can be generalized for the $g$-{non-linear projections}.

\subsection{$\mathcal{C}$-convex bodies and their companions}

First, we extend the Definitions~\ref{D-supp} and \ref{D-visualh}.

\begin{defi}\label{D-kball-sup}\rm
A ball $B^{k+1}(x,r)$ with the boundary sphere $S^{k}(x,r)$ in ${\Bbb{E}}^n$
will be called \textbf{supporting} a body $K\subset{\Bbb{E}}^n$
if $S^{k}(x,r)\cap K\ne\emptyset$ and $B\/^{k+1}(x,r)$ is disjoint from int\,$K$.
 \end{defi}

\begin{defi}\label{D-CVH}\rm
Let $\mathcal{P}_{k}$ be a non-empty family of circular projections onto punctured $(n-k)$-planes~in~${\Bbb{E}}^n$.
We say that subsets $A,B$ of ${\Bbb{E}}^n$ are $\mathcal{P}_k$-$\mathcal{C}$-\textbf{equivalent},
written $A\overset{\mathcal{P}_k}\sim B$, if and only if $f(A)=f(B)$ for all $f\in\mathcal{P}_k$.
 For each subset $A$ of ${\Bbb{E}}^n$ we define $\mathcal{P}_k\langle A\rangle$ to be the union of all subsets $B$ satisfying $B\overset{\mathcal{P}_k}\sim A$.

The $\mathcal{C}$-\textbf{visual hull} $\mathcal{CP}_{k}\langle W\rangle$ of a set $W$ in ${\Bbb{E}}^n$
is the largest set~$V$ such that $f(V)\subseteq f(W)$ for all $f\in\mathcal{P}_{k}$.
 A body $K\subset{\Bbb{E}}^n$ is called $\mathcal{P}_k$-$\mathcal{C}$-\textbf{visual} if $\mathcal{CP}_{k}\langle K\rangle=K$.
\end{defi}

\begin{prop} $\,\mathcal{CP}_{k}\langle K\rangle= K$ if and only if for every $x\notin K$ there is $f_{P(C)}\,{\in}\,\mathcal{P}_{k}$ such that $S^\perp(P(C),x)$ does not intersect~$K$.
\end{prop}

\textbf{Proof} almost completely repeats the proof of Proposition~\ref{P-003}.

1. Let $x\in\mathcal{CP}_k\langle K\rangle\setminus K$. Thus for all circular projections
$f_{P(C)}\in\mathcal{P}_k$ the image of this point, $f_{P(C)}(x)$, belongs to $f_{P(C)}(K)$,
but then the corresponding circles $S^\perp(P(C),x)$ will intersect $K$.

2. Let for any point $x\not\in K$ there exists a circular projection $f_{P(C)}\in\mathcal{P}_k$
such that $S^\perp(P(C),x)$ does not intersect $K$. Thus $f_{P(C)}\not\in\mathcal{P}_k$,
and we may conclude that any point outside of $K$ does not belong to $\mathcal{CP}_k\langle K\rangle$.$\,\square$

\vskip1mm
Consider several useful families $\mathcal{P}_k$.

\begin{exam}\rm Given a set $M\subset{\Bbb{E}}^n$, denote by $\mathcal{P}_{k,M}$ the family of circular projections onto punctured $(n-k)$-planes $P(C)$ in ${\Bbb{E}}^n$ such that $C\in M$.

If a body $K$ belongs to a ball $B(x,r)$, one may consider $M=S^{n-1}(x,r)$ -- a sphere.

The following two examples are used in our study:

1. Given a body $K\subset{\Bbb{E}}^n$, denote by $\mathcal{P}_{k,K,\,\varepsilon}$ a family of circular projections onto punctured $(n{-}k)$-planes $P(C)$ such that $P^\perp(C)\cap K=\emptyset$, and $dist(C, K)\le 1/\varepsilon$.
If we take $M=\{C\in{\Bbb{E}}^n: {\rm dist}\,(C,K)\le 1/\varepsilon\}$ then we obtain the classes of Definition~\ref{D-Kkeps1} in what follows.

2. Given a $k$-plane $L^k\subset{\Bbb{E}}^n$, denote by $\mathcal{P}_{L^k}$ a family of circular projections onto punctured $(n-k)$-planes $P(C)$ in ${\Bbb{E}}^n$ such that $C\in L^k$ and $P(C)\perp L^k$. For $\mathcal{P}_{k}=\mathcal{P}_{L^k}$ in Definition~\ref{D-Kkeps1}
we obtain 4 classes of bodies, that will be denoted by $\mathcal{KP}_{k,i}(L^k),\mathcal{VP}_{k,i}(L^k)$, $i=1,2$.
\end{exam}

In order to study the \textbf{question} (see Introduction), we will extend the Definition~\ref{D-Kki}.

\begin{defi}\label{D-Kkeps1}\rm
A body $K\subset{\Bbb{E}}^n$ is called $\mathcal{C}$-\textbf{convex} of a class ${\mathcal{K}}^{k,\,\varepsilon}_i$ (for some $\varepsilon>0$ and $0\le k\le n-1$) if

\vskip.1mm\hskip-2mm
${\mathcal{K}}^{k,\,\varepsilon}_1$:
{\it for any point $x\in\partial K$ there is a $f_{P(C)}\in\mathcal{P}_{k,K,\,\varepsilon}$ such that
$x\in P(C)$, $|Cx|\,{=}\,1/\varepsilon$
and $f_{P(C)}(x)\in\partial(f_{P(C)}(K))$},

\vskip.1mm\hskip-2mm
${\mathcal{K}}^{k,\,\varepsilon}_2$:
{\it for any point $x\notin K$, dist$(x,K)\le1/\varepsilon$,  there is a $f_{P(C)}\in\mathcal{P}_{k,K,\,\varepsilon}$ such that $x\in P(C)$, $|Cx|\,{=}\,1/\varepsilon$
and $f_{P(C)}(x)$ is disjoint from $f_{P(C)}(K)$}.

\noindent
A body $K\subset{\Bbb{E}}^n$ is called $\mathcal{C}$-\textbf{visible} of a class ${\mathcal{V}}^{\,k,\,\varepsilon}_i$ if

\hskip-2mm
${\mathcal{V}}^{\,k,\,\varepsilon}_1$:
{\it for any ball $\tilde B^{k-m+1}(C,1/\varepsilon)$ supporting $K$ there is a $f_{P(C)}\in\mathcal{P}_{k,K,\,\varepsilon}$
such that $\partial\tilde B\,{\perp}\,S(P(C),1/\varepsilon)$ and $f_{P(C)}(\tilde B)$ is supporting $f_{P(C)}(K)$},

\hskip-2mm
${\mathcal{V}}^{\,k,\,\varepsilon}_2$:\,{\it for any sphere $\tilde S^{k-m}(C,1/\varepsilon)$ disjoint from $K$ and dist\,$(C,K)\le 1/\varepsilon$, there is a $f_{P(C)}\,{\in}\,\mathcal{P}_{k,K,\,\varepsilon}$
 such that $\tilde S\,{\perp}\,S(P(C),1/\varepsilon)$ and $f_{P(C)}(\tilde S)$ is disjoint from $f_{P(C)}(K)$}.

\vskip.5mm\noindent
Here, as above, $0<m\le k<n$, and the $0$-dimensional sphere is a point.
\end{defi}

Similarly one may define the classes $\mathcal{KP}_{k,i}\ (i=1,2)$ of $\mathcal{C}$-\textbf{convex} and $\mathcal{C}$-\textbf{visible} bodies in ${\Bbb{E}}^n$ relative to \textbf{any} family $\mathcal{P}_{k}$ (in particular, for $\mathcal{P}_{L^k}$ and $\mathcal{P}_{k,M}$) of circular projections onto punctured $(n{-}k)$-planes $P(C)$.

One may verify (applying just the set theory arguments) that if a body $K$ is the intersection of connected bodies of a class $\mathcal{K}^{k,\varepsilon}_2$ then $K$ also belongs to $\mathcal{K}^{k,\varepsilon}_2$. This intersection can be disconnected.

\begin{prop}\label{P-separatekeps}
The conditions $K\in\mathcal{K}^{k,\,\varepsilon}_i\ (i=1,2)$ are equivalent to the following ones:

$(i=1)$ \textit{any point $x\in\partial K$ belongs to a
ball $B^{k+1}(C,1/\varepsilon)$ supporting~$K$},

 $(i=2)$ \textit{any point $x\notin K$, dist$(x,K)\le1/\varepsilon$, belongs to
 the boundary  sphere of a ball $B^{k+1}(C,1/\varepsilon)$ disjoint from $K$}.

\noindent
The condition $K\in\mathcal{V}^{k,\varepsilon}_i\,(i=1,2)$ is equivalent to the following ones:

 $(i=1)$
\textit{any ball $B^{k-m+1}(C,1/\varepsilon)$ supporting $K$ belongs to a ball $B^{k+1}(C,1/\varepsilon)$ supporting $K$ as well},

 $(i=2)$ \textit{any sphere $S^{k-m}(C,1/\varepsilon)$ disjoint from $K$
 and dist\,$(C,K)\le 1/\varepsilon$ belongs to the boundary sphere of
 a ball $B^{k+1}(C,1/\varepsilon)$ disjoint from $K$ as well}.
\end{prop}

The proof of Proposition~\ref{P-separatekeps} obviously follows from the definitions.

\begin{rem}\rm
We have introduced in \cite{vr} a "$k=n-1$ version" of the $\mathcal{C}$-convex objects.
 A body $K\subset M$ was called $\varepsilon$-\textbf{convex} of a class $\mathcal{K}^\varepsilon_i$ (for some $\varepsilon>0$)~if

\hskip-1mm
$\mathcal{K}^\varepsilon_1$:
\textit{any point $x\in\partial K$ belongs to a supporting ball (of $K$)
of radius~$1/\varepsilon$},

\hskip-1mm
$\mathcal{K}^\varepsilon_2$:
\textit{any point $x\notin K$ belongs to a ball $B$ of radius $1/\varepsilon$ such that
int\/$K\cap B{=}\emptyset$}.

\noindent
Here only the class $\mathcal{K}^{\varepsilon}_1 $ coincides with $\mathcal{K}^{n-1,\,\varepsilon}_1$.
Simple examples show that corresponding classes $\mathcal{K}^{n-1,\,\varepsilon}_2$ and $\mathcal{K}^{\varepsilon}_2$
are different.
\end{rem}

\begin{theo}\label{P-separate3} The strong inclusions hold for all $\varepsilon>0,\ k<n,\ i=1,2$:
\vskip.5mm
$
 \mathcal{K}^{k,\,\varepsilon}_2\subset\mathcal{K}^{k,\,\varepsilon}_1$,\ \
$
 \mathcal{V}^{k,\,\varepsilon}_2\subset\mathcal{V}^{k,\,\varepsilon}_1$,\ \
and\ \
$\mathcal{V}^{k,\,\varepsilon}_i\subset\mathcal{K}^{k,\,\varepsilon}_i$.
\end{theo}

\textbf{Proof} is similar to the proof of Theorem~1 in \cite{vr}.
\newline
$\mathcal{K}^{k,\,\varepsilon}_2\subset\mathcal{K}^{k,\,\varepsilon}_1$.
Let $ x_* \in \partial K $, $K\in\mathcal{K}^{k,\,\varepsilon}_2$.
Consider a sequence of points $x_i\notin K$ converging to $x_*$.
By condition, for any $i$ there is a sphere $S^k(C_i, 1/\varepsilon)$, containing $x_i$ and disjoint from $K$, moreover, dist$(C_i,K)\le 1/\varepsilon$.
The bounded sequence of points $\{C_i\}$ has a limit point $C_*$, moreover,
dist$(C_*,K)\le 1/\varepsilon$ and
$$
 1/\varepsilon\geq\lim\limits_{i\to\infty}\rho(x_i, C_i)=\rho(x_*,C_*).
$$
Hence, $ x_*\in S^k(C_*, 1/\varepsilon)\,$ and a sphere
$\,S^k (C_*, 1/\varepsilon) $ is disjoint from int\,$K$ and determines a circular projection from a class $\mathcal{K}^{k,\,\varepsilon}_1$.

Other inclusions can be proved analogously.$\,\square$

\begin{figure}
\begin{center}
\includegraphics[scale= 0.3,angle=0]{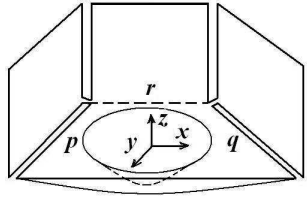}
\caption{\small $K\in\mathcal{K}^{k,\varepsilon}_1\setminus\mathcal{K}^{k,\varepsilon}_2\ (k=1)$.
}
\label{ex-K1234eps}
\end{center}
\end{figure}

\begin{exam}\label{Ex-1-2eps}\rm
One can produce examples of bodies $K\in\mathcal{K}^{k,\,\varepsilon}_1\setminus\mathcal{K}^{k,\,\varepsilon}_2$
applying small deformations to corresponding examples of $K\in\mathcal{K}^{k}_1\setminus\mathcal{K}^{k}_2$,
see Fig.~\ref{ex-K12}, where a circular domain about the coordinate center is replaced by a spherical shell of radius $1/\varepsilon$, see Fig.~\ref{ex-K1234eps}.
\end{exam}

The following proposition was proved in \cite{vr} for $k=n-1$.

\begin{prop}\label{P-03}
If a connected body $K\in\mathcal{K}^{k,\,\varepsilon}_1$ has {\rm diam}\,$K<1/\varepsilon$
then the boundary $\partial K$ is connected as well.
\end{prop}

\textbf{Proof}. Let $K$ be connected and its boundary $\partial K$ is not, then one of components of $\partial K$, denote it by $\partial_1 K$, separates $K$ from "infinity"
(or from ${\Bbb{E}}^n\setminus conv K_\varepsilon$).
Now, $\partial_1 K$ is a boundary of some body $L$ such that diam\,$L={\rm diam\,}(\partial K)={\rm diam\,}K$.
Let $\partial_2 K\neq\partial_1 K$ be another component of $\partial K$ and let $x_1\in\partial_2 K$.
If $K\in\mathcal{K}^{k,\varepsilon}_1$ then $x_1$ should be contained in $L$ and in some supporting ball  $B^{k+1}(C,1/\varepsilon)$, moreover, dist$(C,K)\le 1/\varepsilon$. So, $B^{k+1}(C,1/\varepsilon)\subset L$ and diam\,$L\ge1/\varepsilon$, a contradiction.$\,\square$

\vskip1mm
Next we consider a special case of Definition~\ref{D-CVH}.

\begin{defi}\label{D-epsVH}\rm
The $\varepsilon$-\textbf{visual hull} $\mathcal{CP}_{k,\,\varepsilon}\langle V\rangle$ of a set $V$ in ${\Bbb{E}}^n$
(defined by the family $\mathcal{P}_{k,V,\,\varepsilon}$ of circular projections) is the largest set~$K$ such
that $f(K)\subseteq f(V)$ for all $f\in\mathcal{P}_{k,V,\,\varepsilon}$.
A~body $K\subset{\Bbb{E}}^n$ is called $\varepsilon$-\textbf{visual} (relative to $\mathcal{P}_{k,K,\,\varepsilon}$)
if $\mathcal{CP}_{k,\varepsilon}\langle K\rangle=K$.
\end{defi}

\begin{prop} If diam\,$K<1/\varepsilon$, then $K\in\mathcal{K}^{k,\,\varepsilon}_2$ if and only if $\,\mathcal{CP}_{k,\,\varepsilon}\langle K\rangle=K$.
\end{prop}

\textbf{Proof}. From definition of the class ${\mathcal{K}}^{k,\,\varepsilon}_2$ it follows that for any point $x\notin K$ there is a circular projection $f_{P(C)}\in\mathcal{P}_{k,K,\,\varepsilon}$ such that
$|Cx|\,{=}\,1/\varepsilon$ and $f_{P(C)}(x)$ is disjoint from $f_{P(C)}(K)$. Hence, any point outside of $K$ does not belong to its $\varepsilon$-visual hull $\mathcal{CP}_{k,\,\varepsilon}\langle K\rangle$.$\,\square$

\begin{rem}\rm
For the "large" bodies (diam\,$K\ge1/\varepsilon$) this proposition is wrong:

Let $K = B(0, 3/\varepsilon)\setminus B(0, 2/\varepsilon)$. For any point $ C\in B(0, 2/\varepsilon)$ and for any plane $P(C)$ of any dimension the orthogonal complement $P^{\bot}(C)$ intersects this body $K$ which obviously belongs to ${\mathcal K}^{k,\varepsilon}_2$, for any $k$, $1<k<n$ (for "small" bodies, i.e., $diam(K)<1/\varepsilon$, this claim fails). Hence, for determination of the $\varepsilon$-visual hull ${\mathcal CP}_{k,\varepsilon}(K)$ one should take the centers $C$ outside of the "large" ball $B(0,3/\varepsilon)$, i.e., we have here ${\mathcal CP}_{k,\varepsilon}(K)=B(0, 3/\varepsilon)$.
\end{rem}

\subsection{Reconstruction of $\mathcal{C}$-convex bodies by their circular projections}

\begin{theo}\label{T-003}
If the circular projections of $ V_1, V_2\in{\cal K}^{n-1,\varepsilon}_2$  onto any punctured line $P^1(C)$ coincide, and $diam V_1,\,diam V_2 < 1/(2\varepsilon)$, then $V_1 = V_2$. Here $ C $ is a center of any ball of radius $1/\varepsilon$ supporting $V_1$ or $V_2$.
\end{theo}

\textbf{Proof}.
1. It follows from the conditions that if $ B(C,1/\varepsilon) $ is supporting $V_1$, then it is supporting to $V_2$ and vice versa. In the same way one can verify that $ dist(V_1, V_2) < 1/(2\varepsilon) $.

2. Hence $ conv(V_1) = conv(V_2)$. We denote it by $K$, its diameter is less $ 1/(2\varepsilon) $ as well and any center of supporting ball with radius $ 1/\varepsilon$ can not be contained in $K$.

3. Now, if $ y_2 \in V_2 \subset K$ and $ y_2\notin V_1$, then we
can assume that $ y_2\in int V_2 $;
 $ y_2\in B(C(y_2), 1/\varepsilon)\equiv B_2 $, $B_2$ is disjoint from $ int V_1 $ by definition of the
class $\mathcal{K}_2^{n-1,\varepsilon} $. This ball $B_2$ can not be supporting $ V_1 $, otherwise it should be supporting to $V_2 $. So, $ B_2 $ is disjoint from $V_1$.

4. Let $E(y_2)$ be the equatorial hyperplane of $ B_2 $ orthogonal to the line $ [C(y_2), y_2] $. This hyperplane divides the space ${\Bbb{E}}^n $. Since the bodies $ V_1, V_2 $ are "small", one of these  half-spaces ${\Bbb{E}}^n_+$, contains $ V_1 $ and $ V_2 $, another one ${\Bbb{E}}^n_-$ does not contain them.

5. We shall move $ B_2 $ in the direction $\overrightarrow{y_2,C(y_2)}$ -- outside of $K$ and the result will be denoted by $B_2(t)$. During this displacement the balls $ B_2(t) $  will be contained in the union $ B_2 \cup{\Bbb{E}}^n_-$, so they will never intersect $V_1$. At some moment $ t_1 $ the ball $B_2(t_1)$ will be supporting $ V_2 $. But when $B_2 $ moves outside of $ K $ it can not touch $ V_1 $. This contradiction follows from assumption 3.

Hence $ V_1 = V_2 $.$\,\square$

\begin{rem}\rm For "large" bodies $V_1, V_2$, $diam V_1, diam V_2>1/\varepsilon$ and for positive circular projection Theorem~\ref{T-003} is wrong:

Let $ V_1 \subset{\Bbb{E}}^2 $ be a ring $ 1+2/\varepsilon\leq |x|\leq 3+2/\varepsilon $ and
$ V_2 = V_1 \cup B(0,1) $, where $ B(0,1) $ is a unit ball centered
in the origin. All the balls $ B(Z, 1/\varepsilon) $ centered in $ Z $,
$ |Z| = 1+ 1/\varepsilon$ are supporting both for $ V_1 $ and $ V_2 $ (and for $ B(0,1)$).
Clearly, for all such $Z$ we have $|f_Z|(V_1 ) = |f_Z|(V_2) $.

Theorem~\ref{T-003} also fails for a wider class ${\cal K}^{n-1,\varepsilon}_1$ that can be seen from examples similar to above one.
\end{rem}

\begin{prop}
If a body $ V_1\subset{\Bbb{E}}^n\setminus L^k$
belongs to a class $\mathcal{KP}_{k,2}(L^k)\ (1\le k\le n-1)$,
and $V_2$ is obtained from $V_1$ by some rotation of ${\Bbb{E}}^n$ about $L^k$ (i.e., by some transformation in $SO(n-k)$),
then for any $C\in L^k $ the circular projections $f_{P(C)}$ of the bodies $V_1$ and $V_2$ onto the plane
$P^{n-k}(C)$ are congruent with respect to some rotation of this plane (or are $SO(n-k)$-equivalent).
\end{prop}

\textbf{Proof}. The proof immediately follows from the fact that for any $ C\in L^k $  the sphere $ S^{n-1}(C,r)$ and decomposition of the space $ E^n = L^k\oplus P^{n-k}(C) $ are invariant with respect to
the action of $SO(n-k)$.$\,\square$

\vskip1mm
We will say that a body $K\subset{\Bbb{E}}^n$ \textbf{surrounds} a plane $P^{n-2}$ if
any half-plane $H^{n-1}$ bounded by $P^{n-2}$ intersects $K$.

Here is an example of an answer to the \textbf{question} from the Introduction.

\begin{theo}\label{T-03}
If the bodies $V_1, V_2\subset{\Bbb{E}}^n\setminus L^{n-2}$
belong to a class $\mathcal{KP}_{n-2,2}(L^{n-2})$
and their circular projections onto all planes $P^{2}(C)\ (C\in L^{n-2})$ orthogonal to $L^{n-2}$ are $SO(2)$-equivalent, then $V_1$ can be obtained from $V_2$ by some rotation about $L^{n-2}$.
\end{theo}

\textbf{Proof of Theorem~\ref{T-03}}.
For visuality, we divide the proof into two steps: 3-dimensional and multidimensional.
Let $n-2=1$.

1. Fix some $C_0\in L^1$ and rotate $V_1$ about the axis $L^1$ so that the image $V_1'$ and the body $V_2$ have equal circular projection onto $P^2(C_0)$, i.e., $f_{C_0}(V_1')=f_{C_0}(V_2)$.
These projections are contained in some angle $A_0 C_0 B_0 $ in the plane $ P(C_0)$ and touch the edges of this angle. Hence, the bodies $V_1'$ and $V_2$ belong to the dihedral angle $A_0 C_0 B_0$ with the edge $L^1$ and touch the faces of this angle.

2. If for some $C_1\in L^1$ the circular projections $f_{C_1}(V_1')$ and $f_{C_1}(V_2)$ do not coincide and are equivalent with respect to rotation in the plane $P(C_1 $ by some angle $\phi(C_1)$, $0<\phi(C_1)<2\pi$, then $f_{C_1}(V_1')$ is contained in some plane angle $A_1 C_1 B_1$ and touches its edges, and the circular projection $ f_{C_1}(V_2) $ is contained in some equivalent angle $ A_1' C_1 B_1' $ and touches its edges. Hence, the bodies $V_1'$ and $V_2$ are contained in dihedral angles $A_1 C_1 B_1$, $A_1' C_1 B_1'$, respectively; these angles have common edge $L^1$ and their faces touch these bodies $V_1'$ and $V_2$. This contradicts to the existence of the dihedral angle $A_0 C_0 B_0$ constructed above.

Now let $n-2>1$.

3. Fix as above any $C_0\in L^{n-2}$ and rotate the body $V_1$ about the axis $L^{n-2}$ so that the image $V_1'$ and the body $V_2$ have equal circular projection $f_{C_0}(V_1')=f_{C_0}(V_2)$ in $P^{2}(C_0)$.

4. If for some $C_1\in L^{n-2}$ the circular projections $ f_{C_1}(V_1')$ and $f_{C_1}(V_2)$ do not coincide, but are equivalent with respect to some rotation $\varphi(C_1)\in SO(2)$
of the plane $P^2(C_1)$, then since the bodies $V_1,V_2$ do not surround $L^{n-2}$, their circular projections $V_1'(P^2(C_1))$ and $V_2(P^2(C_1))$ as in 3-dimensional case, belong to
equal plane angle $A_1 C_1 B_1$ and $A_1' C_1 B_1'$ with a common $(n-2)$-dimensional edge
$L^{n-2}$ and touche their edges.
 As above, the bodies $V_1'$, $V_2$ are contained in corresponding dihedral angles $A_1 C_1 B_1$ and $A_1' C_1 B_1'$ with common $(n-2)$-dimensional "edge" $L^{n-2}$ and touch their faces.

These angles are congruent with respect to the rotation by the angle $\varphi(C_1)$ in the plane $P^2(C_1)$. On the other hand the intersections of the faces with the plane $P^{n-k}(C_0)$ generate in that plane the dihedral angles with common vertex $C_0$.
These dihedral angles contain the circular projections $f_{C_0}(V_1') = f_{C_0}(V_2)$ and they are congruent with respect to the same nontrivial rotation by the angle $\varphi(C_1)$, a~contradiction to non-triviality of this angle.$\,\square$

\begin{rem}\rm
(a) The case of $SO(2)$-equivalent orthogonal projections on all 2-planes, where under additional assumption that the projections of bodies have no $SO(2)$-symmetries, it was proven that the bodies $V_1,V_2$ are equal with respect either translation or central symmetry,
was studied in \cite{vgol}.
In contrast to this study, the conditions of Theorem~\ref{T-03} do not contain asymmetry assumptions, but

-- only orthogonal to $L^{n-2}$ planes of circular projection are considered,

-- the centers of these circular projections belong to $L^{n-2}$.

(b) In the same way Theorem~\ref{T-03} can be reformulated for the cases of complex and quaternion Euclidean spaces where one should consider $U(n-k),\,SU(n-k)$ and $Sp(n-k)$-transformations.
In \cite{vgol99} corresponding results were obtained for the case of orthogonal projections onto $(n-k)$-planes of complex Euclidean space for arbitrary $k$ and under assumption that
these projections have no $U(n-k)$ or $SU(n-k)$ symmetries.
\end{rem}



\baselineskip=9pt

\begin{thebibliography}{9999}

\bibitem{resh56}
 Yu.\,G. Reshetnyak,
 {On one generalization of convex surfaces}. Math. Sbornik, 1956, 40 (3), 381--398 (Russian).

\bibitem{vr}
 V.\,P. Golubyatnikov, V.\,Y Rovenski,
{Notes on $\varepsilon$-convex bodies}, submitted to {Israel J. of Math.} (2008), 1--18.

\bibitem{vgol}
 V.\,P. Golubyatnikov,
 {Uniqueness Questions in Reconstruction of Multidimensional Objects from Tomography-Type Projection Data}. VSP BV, The Netherlands, 2000.

\bibitem{MU67}
 G.\,H. Meisters, S.\,M. Ulam,
 {On visual hull of sets}. Proc. Nat. Acad. Sci. USA, 1967, 57, 1172--1174.

\bibitem{LM70}
 D.\,G. Larman, P. Mani,
 {On visual hulls}. Pacific J. Math., 1970, 32(1), 157--171.

 \bibitem{Shefel'_1969}
 S. Shefel',
 {About two classes of $k$-dimensional submanifolds in $n$-dimen\-si\-onal Euclidean space}.
 {Siberian Math. J.}, 1969, 10(2): 459--467.

\bibitem{vgol92}
 V.\,P. Golubyatnikov,
 {Stability problems in certain inverse problems of reconstruction of convex compacta from their projections (English)}. {Siberian Math. J.} 33\,(3), (1992), 409--415.

\bibitem{vgol99}
 V.\,P. Golubyatnikov,
 {On the unique determination of compact convex sets from their projections. The complex case}. {Siberian Math. J.} 40\,(4), (1999), 678--681.

\bibitem{Rov_2006}
 V.\,Yu. Rovenski, {On $(k,\varepsilon)$-saddle submanifolds of Riemannian manifolds}.
 {Geometriae Dedicata}, 121 (2006), 187--203.

\bibitem{gro}
 M. Gromov,
 {Spaces and questions}. Geom. Funct. Anal. (GAFA), Special Volume, Part I (2000), 118--161.

\bibitem{rgar06}
 R. Gardner,
 {Geometric Tomography}, Cambridge University Press, 2006.

\bibitem{Bor82}
 A.\,A. Borisenko,
 {On the extrinsic geometric properties of parabolic surfaces and topological  properties of saddle surfaces in symmetric spaces of rank one}.
 {Math. USSR Sb.}, 42: 297--310, 1982.

\bibitem{LF07}
 S. Lazebnik, Y. Furukawa and J. Ponce,
 {Projective visual hulls}. Int. J. of Computer Vision, 2007, 74\,(2), 137--165.


\end{thebibliography}
\end{document}